\documentclass{amsart}
\usepackage{latexsym}

\usepackage{amsfonts,amssymb,amsmath}
\usepackage{amsthm}
\usepackage{graphicx,color}
\usepackage{color}
\usepackage{amsmath}
\newtheorem{theorem}{Theorem}
\newtheorem{lemma}{Lemma}[section]
\newtheorem{corollary}{Corollary}
\newtheorem{proposition}{Proposition}
\newtheorem{conjecture}{Conjecture}
\newtheorem{definition}{Definition}

\newfont{\bg}{cmr9 scaled\magstep4}
\newcommand{\bigzerol}{\smash{\lower1.0ex\hbox{\bg 0}}}

\newcommand{\R}{\mathbb{R}}
\title{
Generalized distance-squared mappings \\
of the plane into the plane
}
\author{S.~Ichiki}
\address{Dai Nippon Printing Co., Ltd.,
Tokyo 162-8001, Japan}
\email{ichiki-shunsuke-jb@ynu.jp}
\author{T.~Nishimura
}
\address{Research Group of Mathematical Sciences,
Research Institute of Environment and Information Sciences, Yokohama
National University, Yokohama 240-8501, Japan}
\email{nishimura-takashi-yx@ynu.jp}
\author{R.~Oset~Sinha}
\address
{Department of Mathematics,   Federal University of Sao Carlos, Sao
Carlos, SP, Brazil} \email{Raul.Oset@uv.es}
\author{M.~A.~S.~Ruas}
\address{ICMC,
University of Sao Paulo, Sao Carlos, SP, Brazil}
\email{maasruas@icmc.usp.br}
\begin{document}
\date{}
\begin{abstract}
We define generalized distance-squared mappings, and we concentrate
on the plane to plane case. We classify generalized distance-squared
mappings of the plane into the plane in a recognizable way.
\end{abstract}
\maketitle \noindent
\section{Introduction}\label{Introduction}
Let $n$ be a positive integer.   For any positive integer $k$, let
$p_0, p_1, \ldots, p_k$ be $(k+1)$ points of $\mathbb{R}^{n+1}$. Set
$p_i=(p_{i0}, p_{i1}, \ldots, p_{in})$ $(0\le i\le k)$. Let
$A=(a_{ij})_{0\le i\le k, 0\le j\le n}$ be a $(k+1)\times (n+1)$
matrix with non-zero entries. Then, the  following mapping $G_{(p_0,
p_1, \ldots, p_k, A)}: \mathbb{R}^{n+1}\to \mathbb{R}^{k+1}$ is
called a {\it generalized distance-squared mapping}:
\[
G_{(p_0, p_1, \ldots, p_k, A)}(x)=\left( \sum_{j=0}^n
a_{0j}(x_j-p_{0j})^2, \sum_{j=0}^n a_{1j}(x_j-p_{1j})^2, \ldots,
\sum_{j=0}^n a_{nj}(x_j-p_{nj})^2 \right),
\]
where $x=(x_0, x_1, \ldots, x_n)$. In \cite{ichikinishimura} (resp.,
\cite{ichikinishimura2}), the first and the second authors gave a
classification result on distance-squared mapping $D_{(p_0, p_1,
\ldots, p_k)}$ (resp., Lorentzian distance-squared mappings
$L_{(p_0, p_1, \ldots, p_k)}$) which is the mapping $G_{(p_0, p_1,
\ldots, p_k,A)}$ in the case that each entry of $A$ is $1$ (resp.,
in the case of $a_{i0}=-1 $ and $a_{ij}=1$ if $j\ne 0$). In these
cases, the rank of $A$ is $1$. Generalized distance-squared mappings
are a useful tool in the applications of singularity theory to
differential geometry. Their singularities give information on the
contacts amongst the families of quadrics defined by the components
of $G_{(p_0,p_1,\ldots,p_k,A)}$. It is therefore natural to classify
maps $G_{(p_0,p_1,\ldots,p_k,A)}$, which is the main theme of this
paper.
\par
Two mappings $f_i: \mathbb{R}^{n+1}\to \mathbb{R}^{k+1}$
$(i=1, 2)$ are said to be $\mathcal{A}$-equivalent if there exist
$C^\infty$ diffeomorphisms $h: \mathbb{R}^{n+1}\to \mathbb{R}^{n+1}$
and $H:\mathbb{R}^{k+1}\to \mathbb{R}^{k+1} $
such that the identity $f_1=H\circ f_2\circ h$ holds.
It turns out that new $\mathcal{A}$-classes occur even in the
simplest equidimensional case (see Theorem \ref{theorem 1}).
\begin{definition}\label{definition 1}
{\rm
\begin{enumerate}
\item
Let $\Phi_{n+1}: \mathbb{R}^{n+1}\to \mathbb{R}^{n+1}$ denote the
following mapping:
\[
\Phi_{n+1}(x_0, x_1, \ldots, x_n) =\left(x_0, x_1, \ldots, x_{n-1},
x_n^2\right).
\]
When a map-germ $f: (\mathbb{R}^{n+1}, q)\to (\mathbb{R}^{n+1},
f(q))$ is $\mathcal{A}$-equivalent to $\Phi_{n+1}:
(\mathbb{R}^{n+1}, 0)\to (\mathbb{R}^{n+1},0)$, the point $q\in
\mathbb{R}^{n+1}$ is said to be a {\it fold point of $f$}.
\item
Let $\Gamma_{n+1}: \mathbb{R}^{n+1}\to \mathbb{R}^{n+1}$ denote the
following mapping:
\[
\Gamma_{n+1}(x_0, x_1, \ldots, x_n) =\left(x_0, x_1, \ldots,
x_{n-1}, x_n^3+x_0x_n\right).
\]
When a map-germ $f: (\mathbb{R}^{n+1}, q)\to (\mathbb{R}^{n+1},
f(q))$ is $\mathcal{A}$-equivalent to $\Gamma_{n+1}:
(\mathbb{R}^{n+1}, 0)\to (\mathbb{R}^{n+1},0)$, the point $q\in
\mathbb{R}^{n+1}$ is said to be a {\it cusp point of $f$}.
\end{enumerate}
}
\end{definition}
\noindent It is known that both $\Phi_{n+1}, \Gamma_{n+1}$ are
proper and stable mappings (for instance see \cite{arnoldetall}).
\par
\begin{figure}[htbp]
\begin{center}
\includegraphics[width=1.1\linewidth]{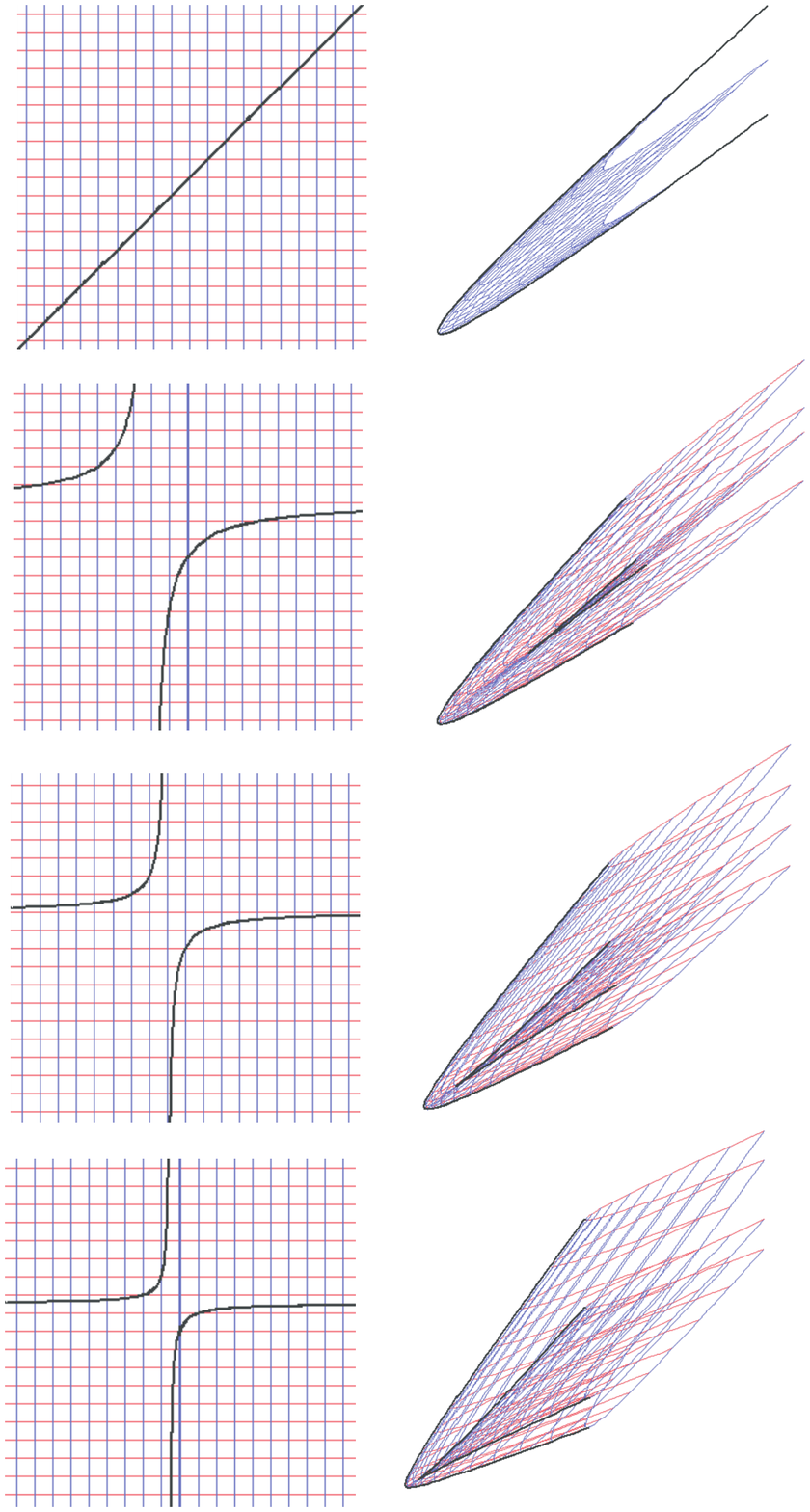}
\caption{Various figures of $S(G_{(p_0, p_1, A_2)})$ and $G_{(p_0,
p_1, A_2)}(S(G_{(p_0,p_1, A_2)}))$ } \label{figure 1}
\end{center}
\end{figure}

In the case of $n=k$, $D_{(p_0, p_1, \ldots, p_n)}$ and $L_{(p_0,
p_1, \ldots, p_n)}$ can be generically characterized as follows.
\begin{proposition}[\cite{ichikinishimura, ichikinishimura2}]\label{proposition 1}
Let $p_0, p_1, \ldots, p_n$ be $(n+1)$-points of $\mathbb{R}^{n+1}$
such that the dimension of
$\sum_{i=1}^n\mathbb{R}\overrightarrow{p_0p_i}$ is $n$. Then, the
following hold:
\begin{enumerate}
\item The distance-squared mapping $D_{(p_0, p_1, \ldots, p_n)}: \mathbb{R}^{n+1}\to \mathbb{R}^{n+1}$
is $\mathcal{A}$-equivalent to $\Phi_{n+1}$.
\item The Lorentzian distance-squared mapping $L_{(p_0, p_1, \ldots, p_n)}: \mathbb{R}^{n+1}\to \mathbb{R}^{n+1}$
is $\mathcal{A}$-equivalent to $\Phi_{n+1}$.
\end{enumerate}
\end{proposition}
For generalized distance-squared mappings,
it is natural to expect that for generic $p_0, p_1, \ldots, p_n$,
$G_{(p_0, p_1, \ldots, p_n, A)}$ is proper and stable, and the rank
of $A$ is a complete invariant of $\mathcal{A}$-types.
\begin{conjecture}\label{conjecture 1}
Let $A_k$
be an $(n+1)\times (n+1)$ matrix of rank $k$ with non-zero entries
$(1\le k\le (n+1))$. Then, there exists a subset $\Sigma\subset
(\mathbb{R}^{n+1})^{n+1}$ of Lebesgue measure zero such that for any
$(p_0, p_1, \ldots, p_n)\in  (\mathbb{R}^{n+1})^{n+1}- \Sigma$, the
following hold:
\begin{enumerate}
\item For any $k$ $(1\le k\le (n+1))$,
the generalized distance-squared mapping $G_{(p_0, p_1, \ldots, p_n,
A_k)}$ is proper and stable.
\item For any two integers $k_1, k_2$ such that $1\le k_1<k_2\le (n+1)$,
$G_{(p_0, p_1, \ldots, p_n, A_{k_2})}$ is not
$\mathcal{A}$-equivalent to $G_{(p_0, p_1, \ldots, p_n, A_{k_1})}$.
\item Let $B_k$ be an $(n+1)\times (n+1)$ matrix of rank $k$ with non-zero entries $(1\le k\le (n+1))$ and
let  $(q_0, q_1, \ldots, q_n)$ be in $ (\mathbb{R}^{n+1})^{n+1}-
\Sigma$.     Then, $G_{(p_0, p_1, \ldots, p_n, A_k)}$ is
$\mathcal{A}$-equivalent to $G_{(q_0, q_1, \ldots, q_n, B_k)}$ for
any $k$.
\end{enumerate}
\end{conjecture}
The main purpose of this paper is to give the affirmative answer to
Conjecture \ref{conjecture 1} in the case $n=1$ as follows.
\begin{theorem} \label{theorem 1}
Let $((x_0, y_0), (x_1, y_1))$ be the standard coordinates of
$(\mathbb{R}^2)^2$ and let $\Sigma$ be the hypersurface in
$(\mathbb{R}^2)^2$ defined by $(x_0-x_1)(y_0-y_1)=0$. Let $(p_0,
p_1)$ be a point in $(\mathbb{R}^2)^2 - \Sigma$ and let $A_k$ be a
$2\times 2$ matrix of rank $k$ with non-zero entries (k=1, 2). Then,
the following hold:
\begin{enumerate}
\item
The mapping $G_{(p_0, p_1, A_1)}$ is $\mathcal{A}$-equivalent to
$\Phi_2$.
\item
The mapping $G_{(p_0, p_1, A_2)}$ is proper and stable, and it is
not $\mathcal{A}$-equivalent to  $G_{(p_0, p_1, A_1)}$.
\item Let $B_2$ be a $2\times 2$ matrix of rank $2$ with non-zero entries and let
$(q_0, q_1)$ be a point in $(\mathbb{R}^2)^2 - \Sigma$. Then,
$G_{(p_0, p_1, A_2)}$ is $\mathcal{A}$-equivalent to $G_{(q_0, q_1,
B_2)}$.
\end{enumerate}
\end{theorem}
\noindent
It turns out that there exists only one cusp point in the singular
set of $G_{(p_0, p_1, A_2)}$ (see Proposition \ref{proposition 2}).
Note that $\Gamma_2 $ is a cubic polynomial. Thus, Theorem
\ref{theorem 1} implies that there exists only one point at which
the $\mathcal{A}$-type of the germ of the quadratic polynomial
$G_{(p_0, p_1, A_2)}$ must be expressed in a  cubic form, although
there are no such points for  $G_{(p_0, p_1, A_1)}$. In Figure
\ref{figure 1} which were created by using \cite{montesinos}, it can
be seen how the singular set of $G_{(p_0, p_1, A)}$ (denoted by
$S(G_{(p_0, p_1, A)})$) and its image $G_{(p_0, p_1,
A_2)}(S(G_{(p_0, p_1, A)}))$ change if the matrix $A$ moves for
fixed $p_0=(0,0), p_1=(1,1)$. The first row of Figure \ref{figure 1}
is for $a_{00}=a_{01}=a_{10}=a_{11}=1$, the second row is for
$a_{00}=a_{01}=a_{10}=1, a_{11}=1.5$,
 the third row is for $a_{00}=a_{01}=a_{10}=1, a_{11}=2$ and the last row is for $a_{00}=a_{01}=a_{10}=1, a_{11}=3$.
The keys for proving Theorem \ref{theorem 1} are the following two
propositions.
\begin{proposition}\label{proposition 2}
Let $A_2$ be a $2\times 2$ matrix of rank two with non-zero entries.
Let $p_0, p_1$ be two points of $\mathbb{R}^2$ satisfying $(p_0,
p_1)\in (\mathbb{R}^2)^2-\Sigma
$, where $\Sigma\subset (\mathbb{R}^2)^2$ is the hypersurface
defined in Theorem \ref{theorem 1}. Then, the following hold:
\begin{enumerate}
\item The singular set $S(G_{(p_0, p_1, A_2)})$ is an equilateral hyperbola.
\item Any point of $S(G_{(p_0, p_1, A_2)})$is a fold point except for one.
\item The exceptional point given in (2) is a cusp point.
\end{enumerate}
\end{proposition}
\begin{proposition}\label{proposition 3}
Let $A_2$ be a $2\times 2$ matrix of rank two with non-zero entries.
Let $p_0, p_1$ be two points of $\mathbb{R}^2$ satisfying $(p_0,
p_1)\in (\mathbb{R}^2)^2-\Sigma
$, where $\Sigma\subset (\mathbb{R}^2)^2$ is the hypersurface
defined in Theorem \ref{theorem 1}. Then, for any positive real
numbers $a, b$ $(a\ne b)$, there exists a point $q=(q_0, q_1)\in
\mathbb{R}^2$ such that $(q, (0,0))\in (\mathbb{R}^2)^2-\Sigma$ and
$G_{(p_0, p_1, A_2)}$ is $\mathcal{A}$-equivalent to $F_q:
\mathbb{R}^2\to \mathbb{R}^2$ defined by
\[
F_q(x,y)=\left((x-q_0)^2+(y-q_1)^2, ax^2+by^2\right).
\]

\end{proposition}
%
\par
\bigskip
Propositions \ref{proposition 2}, \ref{proposition 3} and Theorem
\ref{theorem 1}
are
proved in Sections \ref{section 2}, \ref{section 3} and \ref{section
4}
respectively. Finally, in Section \ref{section 5} we give a
geometric interpretation of the singularities of $G_{(p_0,p_1,A_2)}$
by which it is clearly explained why the cusp point is only one in
Proposition \ref{proposition 2}.
\section{Proof of Proposition \ref{proposition 2}}\label{section 2}
\subsection{Proof of (1) of Proposition \ref{proposition 2}}\label{subsection 2.1}
The Jacobian matrix $JG_{(p_0,p_1, A_2)}(x,y)$ of the mapping
$G_{(p_0,p_1, A_2)}$ at $(x, y)$ is the following:
\begin{eqnarray*}
JG_{(p_0,p_1, A_2)}(x,y)=
\begin{pmatrix}
2a_{00}(x-p_{01}) & 2a_{01}(y-p_{11}) \\
2a_{10}(x-p_{10}) & 2a_{11}(y-p_{11})
\end{pmatrix},
\end{eqnarray*}
and we get
\begin{eqnarray*}
\det JG_{(p_0,p_1, A_2)}(x,y)=
4\left(\left( \left(a_{00}a_{11}-a_{01}a_{10}\right)x+\left(-a_{00}a_{11}p_{01}+a_{01}a_{10}p_{10}\right)\right)y\right.\\
\left.+\left(-a_{00}a_{11}p_{11}+a_{01}a_{10}p_{01}\right)x+\left(a_{00}a_{11}p_{00}p_{11}-a_{01}a_{10}p_{01}p_{10}\right)\right).
\end{eqnarray*}
Hence, the Jacobian matrix of the function $\det JG_{(p_0,p_1,
A_2)}:\R^2\to\R$ is the following:
\begin{eqnarray*}
J(\det JG_{(p_0,p_1, A_2)})(x,y)=4\left((a_{00}a_{11}-a_{01}a_{10})y+(-a_{00}a_{11}p_{11}+a_{01}a_{10}p_{01})\right.,\\
\left.(a_{00}a_{11}-a_{01}a_{10})x+(-a_{00}a_{11}p_{00}+a_{01}a_{10}p_{10})\right).
\end{eqnarray*}
We show that there exist no points in $\R^2$ such that both $\det
JG_{(p_0,p_1,A_2)}(x,y)=0$ and
$(a_{00}a_{11}-a_{01}a_{10})x+(-a_{00}a_{11}p_{00}+a_{01}a_{10}p_{10})=0$
hold. Suppose that there exists a point
$(\widetilde{x},\widetilde{y})$ in $\R^2$ such that both $\det
JG_{(p_0,p_1, A_2)}(\widetilde{x},\widetilde{y})=0$ and
$(a_{00}a_{11}-a_{01}a_{10})\widetilde{x}+(-a_{00}a_{11}p_{00}+a_{01}a_{10}p_{10})=0$
hold. Then, by
$\widetilde{x}=\frac{a_{00}a_{11}p_{00}-a_{01}a_{10}p_{10}}{a_{00}a_{11}-a_{01}a_{10}}$,
\begin{eqnarray*}
\det JG_{(p_0,p_1, A_2)} (\widetilde{x},\widetilde{y})&=&
\frac{4a_{00}a_{01}a_{10}a_{11}(p_{00}-p_{10})(p_{01}-p_{11})}{a_{00}a_{11}-a_{01}a_{10}}\\
&=&0.
\end{eqnarray*}
This contradicts the assumption that $p_{00}\not=p_{10}$ and
$p_{01}\not=p_{11}$. Hence, if $\det JG_{(p_0,p_1, A_2)}(x,y)=0$,
then it follows that
$(a_{00}a_{11}-a_{01}a_{10})x+(-a_{00}a_{11}p_{00}+a_{01}a_{10}p_{10})\not=0$
and $J(\det JG_{(p_0,p_1, A_2)})(x,y)\not=(0,0)$.
We see that $\det JG_{(p_0,p_1, A_2)}(x,y)=0$ if and only if
\begin{eqnarray*}
((a_{00}a_{11}-a_{01}a_{10})x+(-a_{00}a_{11}p_{00}+a_{01}a_{10}p_{10}))y
+(-a_{00}a_{11}p_{11}+a_{01}a_{10}p_{01})x\\+(a_{00}a_{11}p_{00}p_{11}-a_{01}a_{10}p_{01}p_{10})=0.
\end{eqnarray*}
Then, since
$(a_{00}a_{11}-a_{01}a_{10})x+(-a_{00}a_{11}p_{00}+a_{01}a_{10}p_{10})\not=0$,
it follows that
\begin{eqnarray*}
y=-\frac{(-a_{00}a_{11}p_{11}+a_{01}a_{10}p_{01})x+(a_{00}a_{11}p_{00}p_{11}-a_{01}a_{10}p_{01}p_{10})}
{(a_{00}a_{11}-a_{01}a_{10})x+(-a_{00}a_{11}p_{00}+a_{01}a_{10}p_{10})}.
\end{eqnarray*}
Note that the asymptote of the hyperbola is the following:
\begin{eqnarray*}
x=\frac{a_{00}a_{11}p_{00}-a_{01}a_{10}p_{10}}{a_{00}a_{11}-a_{01}a_{10}},\
y=\frac{a_{00}a_{11}p_{11}-a_{01}a_{10}p_{01}}{a_{00}a_{11}-a_{01}a_{10}}.
\end{eqnarray*}
Thus, the singular set $S(G_{(p_0,p_1, A_2)})$ is an equilateral
hyperbola:
\begin{eqnarray*}
S(G_{(p_0,p_1, A_2)})&=&\left\{(x,\phi (x))\in \R^2\right\},
\end{eqnarray*}
where
\begin{eqnarray*}
\phi
(x)&=&-\frac{(-a_{00}a_{11}p_{11}+a_{01}a_{10}p_{01})x+(a_{00}a_{11}p_{00}p_{11}-a_{01}a_{10}p_{01}p_{10})}
{(a_{00}a_{11}-a_{01}a_{10})x+(-a_{00}a_{11}p_{00}+a_{01}a_{10}p_{10})}\\
&=&\frac{-\frac{a_{00}a_{01}a_{10}a_{11}(p_{00}-p_{10})(p_{01}-p_{11})}{(a_{00}a_{11}-a_{01}a_{10})^2}}
{x+\frac{-a_{00}a_{11}p_{00}+a_{01}a_{10}p_{10}}{a_{00}a_{11}-a_{01}a_{10}}}
+\frac{a_{00}a_{11}p_{11}-a_{01}a_{10}p_{01}}{a_{00}a_{11}-a_{01}a_{10}}.
\end{eqnarray*}
\hfill $\Box$
\subsection{Proof of (2) of Proposition \ref{proposition 2}}\label{subsection 2.2}
For the proofs of (2), (3) of Proposition \ref{proposition 2}, the
following lemma is needed.
\begin{lemma}[\cite{whitney}]\label{lemma 1}
Let $f$ be a mapping of the plane into the plane, and let
$q=(q_0,q_1)\in \R^2$ be a singular point of $f$. Suppose that
$J(\det Jf)_q\not=(0,0)$. Put $J(det Jf)_q=(r_0,r_1)$. By a linear
coordinate transformation if necessary, we may assume that
$r_1\not=0$. Then, by the implicit function theorem, there exists a
function $y=\phi (x)$ defined on an open interval containing $q_0$
satisfying $det Jf(x,\phi (x))=0$. Then, $q$ is a fold point of $f$
if and only if
\begin{eqnarray*}
\left.\frac{df(x,\phi (x))}{dx}\right|_{x=q_0}\not=0,
\end{eqnarray*}
and $q$ is a cusp point of $f$ if and only if
\begin{eqnarray*}
\left.\frac{df(x,\phi (x))}{dx}\right|_{x=q_0}=0, \quad
\left.\frac{d^2f(x,\phi (x))}{dx^2}\right|_{x=q_0}\not=0.
\end{eqnarray*}
\end{lemma}
\par
\medskip
We use the notations given in Subsection \ref{subsection 2.1}. By
calculations, $G_{(p_0,p_1, A_2)}(x,\phi (x))$ can be written as
follows:
\begin{eqnarray*}
&&G_{(p_0,p_1,A_2)}(x,\phi (x))=\\
&&\left(a_{00}(x-p_{00})^2+a_{01}(a_{00}a_{11}(p_{01}-p_{11}))^2
\left(\frac{x-p_{00}}{(a_{00}a_{11}-a_{01}a_{10})x+(-a_{00}a_{11}p_{00}+a_{01}a_{10}p_{10})}\right)^2\right.,\\
&&\left.a_{10}(x-p_{10})^2+a_{11}(a_{01}a_{10}(p_{01}-p_{11}))^2
\left(\frac{x-p_{10}}{(a_{00}a_{11}-a_{01}a_{10})x+(-a_{00}a_{11}p_{00}+a_{01}a_{10}p_{10})}\right)^2
\right).
\end{eqnarray*}
Next, by differentiating $G_{(p_0,p_1, A_2)}(x,\phi (x))$, we have
the following:
\begin{eqnarray*}
&&\frac{dG_{(p_0,p_1, A_2)}(x,\phi (x))}{dx}=\\
&&2\left(a_{00}(x-p_{00})+
\frac{a_{01}^2a_{10}(p_{10}-p_{00})(a_{00}a_{11}(p_{01}-p_{11}))^2(x-p_{00})}{((a_{00}a_{11}-a_{01}a_{10})x+(-a_{00}a_{11}p_{00}+a_{01}a_{10}p_{10}))^3}\right.,\\
&&\left.a_{10}(x-p_{10})+
\frac{a_{11}^2a_{00}(p_{10}-p_{00})(a_{01}a_{10}(p_{01}-p_{11}))^2(x-p_{10})}{((a_{00}a_{11}-a_{01}a_{10})x+(-a_{00}a_{11}p_{00}+a_{01}a_{10}p_{10}))^3}\right).
\end{eqnarray*}
Now, we find the solution to the equation $\frac{dG_{(p_0,p_1,
A_2)}(x,\phi (x))}{dx}=0$. We have
\[
\frac{dG_{(p_0,p_1, A_2)}(x,\phi (x))}{dx}=0
\]
if and only if both the following two hold:
\begin{align}
\tag{2.1} a_{00}(x-p_{00})+
\frac{a_{01}^2a_{10}(p_{10}-p_{00})(a_{00}a_{11}(p_{01}-p_{11}))^2(x-p_{00})}
{((a_{00}a_{11}-a_{01}a_{10})x+(-a_{00}a_{11}p_{00}+a_{01}a_{10}p_{10}))^3}&=0,\\
\tag{2.2} a_{10}(x-p_{10})+
\frac{a_{11}^2a_{00}(p_{10}-p_{00})(a_{01}a_{10}(p_{01}-p_{11}))^2(x-p_{10})}
{((a_{00}a_{11}-a_{01}a_{10})x+(-a_{00}a_{11}p_{00}+a_{01}a_{10}p_{10}))^3}&=0.
\end{align}
Furthermore, since $a_{00}\not=0$ and $a_{10}\not=0$, by dividing
$(2.1)$ by $a_{00}$ and dividing $(2.2)$ by $a_{10}$, we get the
following:
\begin{align}
\tag{2.3}
(x-p_{00})\left(1+\frac{a_{00}a_{01}^2a_{10}a_{11}^2(p_{10}-p_{00})(p_{01}-p_{11})^2}
{((a_{00}a_{11}-a_{01}a_{10})x+(-a_{00}a_{11}p_{00}+a_{01}a_{10}p_{10}))^3}\right)&=0,\\
\tag{2.4}
(x-p_{10})\left(1+\frac{a_{00}a_{01}^2a_{10}a_{11}^2(p_{10}-p_{00})(p_{01}-p_{11})^2}
{((a_{00}a_{11}-a_{01}a_{10})x+(-a_{00}a_{11}p_{00}+a_{01}a_{10}p_{10}))^3}\right)&=0.
\end{align}
By $(2.3)$ and $(2.4)$ and the assumption $p_{00}\not=p_{10}$, we
obtain
\begin{eqnarray*}
1+\frac{a_{00}a_{01}^2a_{10}a_{11}^2(p_{10}-p_{00})(p_{01}-p_{11})^2}
{((a_{00}a_{11}-a_{01}a_{10})x+(-a_{00}a_{11}p_{00}+a_{01}a_{10}p_{10}))^3}=0.
\end{eqnarray*}
Thus, the solution to the equation $\frac{dG_{(p_0,p_1, A_2)}(x,\phi
(x))}{dx}=0$ is the following:
\begin{eqnarray*}
x=\frac{(a_{00}a_{01}^2a_{10}a_{11}^2(p_{00}-p_{10})(p_{01}-p_{11})^2)^\frac{1}{3}+(a_{00}a_{11}p_{00}-a_{01}a_{10}p_{10})}
{a_{00}a_{11}-a_{01}a_{10}}.
\end{eqnarray*}
Denote this solution by $\widetilde{x}$.  Then, it is easy to see
that
\begin{eqnarray*}
\phi
(\widetilde{x})=\frac{-(a_{00}^2a_{01}a_{10}^2a_{11}(p_{00}-p_{10})^2(p_{01}-p_{11}))^\frac{1}{3}
+(a_{00}a_{11}p_{11}-a_{01}a_{10}p_{01})}
{a_{00}a_{11}-a_{01}a_{10}}.
\end{eqnarray*}
Therefore, the point satisfying $\frac{dG_{(p_0,p_1, A_2)}(x,\phi
(x))}{dx}=0$ is the only one point $(\widetilde{x},\phi
(\widetilde{x}))$. \hfill $\Box$
\subsection{Proof of (3) of Proposition \ref{proposition 2}}\label{subsection 2.3}
We use the notations given in Subsections \ref{subsection 2.1} and
\ref{subsection 2.2}. By differentiating $\frac{dG_{(p_0,p_1,
A_2)}(x,\phi (x))}{dx}$, we have the following:
\begin{eqnarray*}
&&\frac{d^2G_{(p_0,p_1, A_2)}(x,\phi (x))}{dx^2}=\\
&&2\left(a_{00}+a_{01}^2a_{10}(p_{10}-p_{00})(a_{00}a_{11}(p_{01}-p_{11}))^2f_0(x),\right.\\
&&\left.a_{10}+a_{11}^2a_{00}(p_{10}-p_{00})(a_{01}a_{10}(p_{01}-p_{11}))^2f_1(x)\right),
\end{eqnarray*}
where
\begin{eqnarray*}
f_0(x)&=&\frac{(a_{00}a_{11}-a_{01}a_{10})(-2x+3p_{00})+(-a_{00}a_{11}p_{00}+a_{01}a_{10}p_{10})}
{((a_{00}a_{11}-a_{01}a_{10})x+(-a_{00}a_{11}p_{00}+a_{01}a_{10}p_{10}))^4},\\
f_1(x)&=&\frac{(a_{00}a_{11}-a_{01}a_{10})(-2x+3p_{10})+(-a_{00}a_{11}p_{00}+a_{01}a_{10}p_{10})}
{((a_{00}a_{11}-a_{01}a_{10})x+(-a_{00}a_{11}p_{00}+a_{01}a_{10}p_{10}))^4}
\end{eqnarray*}
Now, we show that there exist no points such that
$\frac{d^2G_{(p_0,p_1, A_2)}(x,\phi (x))}{dx^2}=0$. Suppose that
there exist a point $\widetilde{x}_0$ such that
$\frac{d^2G_{(p_0,p_1, A_2)}(\widetilde{x}_0,\phi
(\widetilde{x}_0))}{dx^2}=0$. Then, we see that
\begin{align}
\tag{2.5} a_{00}+a_{01}^2a_{10}(p_{10}-p_{00})(a_{00}a_{11}(p_{01}-p_{11}))^2f_0(x)&=0,\\
\tag{2.6}
a_{10}+a_{11}^2a_{00}(p_{10}-p_{00})(a_{01}a_{10}(p_{01}-p_{11}))^2f_1(x)&=0.
\end{align}
Since $a_{00}\not=0$ and $a_{10}\not=0$, by dividing $(2.5)$ by
$a_{00}$ and dividing $(2.6)$ by $a_{10}$, we have
\begin{align}
\tag{2.7} 1+a_{00}a_{01}^2a_{10}a_{11}^2(p_{10}-p_{00})(p_{01}-p_{11})^2f_0(x)=0,\\
\tag{2.8}
1+a_{00}a_{01}^2a_{10}a_{11}^2(p_{10}-p_{00})(p_{01}-p_{11})^2f_1(x)=0.
\end{align}
By $(2.7)$ and $(2.8)$, we get
\begin{eqnarray*}
a_{00}a_{01}^2a_{10}a_{11}^2(p_{10}-p_{00})(p_{01}-p_{11})^2(f_0(x)-f_1(x))=0.
\end{eqnarray*}
Since
\begin{eqnarray*}
f_0(x)-f_1(x)=-\frac{3(a_{00}a_{11}-a_{01}a_{10})(p_{10}-p_{00})}
{((a_{00}a_{11}-a_{01}a_{10})x+(-a_{00}a_{11}p_{00}+a_{01}a_{10}p_{10}))^4},
\end{eqnarray*}
we obtain
\begin{eqnarray*}
-\frac{3a_{00}a_{01}^2a_{10}a_{11}^2(a_{00}a_{11}-a_{01}a_{10})(p_{10}-p_{00})^2(p_{01}-p_{11})^2}
{((a_{00}a_{11}-a_{01}a_{10})x+(-a_{00}a_{11}p_{00}+a_{01}a_{10}p_{10}))^4}=0.
\end{eqnarray*}
However, this contradicts the assumption that
$a_{00}a_{01}a_{10}a_{11}(a_{00}a_{11}-a_{01}a_{10})(p_{10}-p_{00})(p_{01}-p_{11})\not=0$.
Therefore, by Lemma \ref{lemma 1}, it follows that $S(G_{(p_0,p_1,
A_2)})$ has only one point $(\widetilde{x},\phi (\widetilde{x}))$ as
a cusp point. \hfill $\Box$
\section{Proof of Proposition \ref{proposition 3}}\label{section 3}
By the parallel translation of the source space $(x,y)\mapsto
(x+p_{00}, y+p_{01})$, $G_{(p_0, p_1, A_2)}$ is
$\mathcal{A}$-equivalent to the following mapping:
\[
(x, y)\mapsto (a_{00}x^2+a_{01}y^2,
a_{10}(x-\widetilde{p}_0)^2+a_{11}(y-\widetilde{p}_1)^2),
\leqno{(3.1)}
\]
where $\widetilde{p}_i=p_{1i}-p_{0i}$ $(0\le i\le 1)$. Note that
$\widetilde{p}_0\widetilde{p}_1\ne 0$. By the parallel translation
of the target space $(X, Y)\mapsto (X-a_{10}\widetilde{p}_0^2,
Y-a_{11}\widetilde{p}_1^2)$, (3.1) is $\mathcal{A}$-equivalent to
the following mapping:
\[
(x, y)\mapsto (a_{00}x^2+a_{01}y^2,
a_{10}x^2-2a_{10}\widetilde{p}_0x+a_{11}y^2-2a_{11}\widetilde{p}_1y).
\leqno{(3.2)}
\]
Set $p(x,y)=a_{00}x^2+a_{01}y^2$,
$q(x,y)=a_{10}x^2-2a_{10}\widetilde{p}_0x+a_{11}y^2-2a_{11}\widetilde{p}_1y$.
Let $M= \left(
\begin{array}{cc}
\alpha_{0} & \alpha_1 \\
\beta_0 & \beta_1
\end{array}\right)
$ denote the $2\times 2$ matrix defined by
\[
M= \left(
\begin{array}{cc}
a_{00} & a_{10} \\
a_{01} & a_{11}
\end{array}
\right)^{-1} \left(
\begin{array}{cc}
1 & a \\
1 & b
\end{array}
\right).
\]
Then $M$ is a regular matrix and we have the following:
\begin{eqnarray*}
{ } & { } & \left(p(x,y), q(x,y)\right)M \\
{} & = &
\left(x^2-2\beta_0a_{10}\widetilde{p}_0x+y^2-2\beta_0a_{11}\widetilde{p}_1y,
ax^2-2\beta_1a_{10}\widetilde{p}_0x+by^2-2\beta_1a_{11}\widetilde{p}_1y\right),
\end{eqnarray*}
which is $\mathcal{A}$-equivalent to
\[
(x,y)\mapsto \left(
(x-\beta_0a_{10}\widetilde{p}_0)^2+(y-\beta_0a_{11}\widetilde{p}_1)^2,
a(x-\frac{1}{a}\frac{1}{2}\beta_1a_{10}\widetilde{p}_0)^2+
b(y-\frac{1}{b}\beta_1a_{11}\widetilde{p}_1)^2 \right).
\]
Therefore, $G_{(p_0, p_1, A_2)}$ is $\mathcal{A}$-equivalent to
$F_{q}$ where $q= \left( \beta_0a_{10}\widetilde{p}_0-
\frac{1}{a}\beta_1a_{10}\widetilde{p}_0,\beta_0a_{11}\widetilde{p}_1-\frac{1}{b}\beta_1a_{11}\widetilde{p}_1
\right)$. Finally, we show that $(q, (0,0))\in
(\mathbb{R}^2)^2-\Sigma$,
\begin{eqnarray*}
\beta_0a_{10}\widetilde{p}_0-
\frac{1}{a}\beta_1a_{10}\widetilde{p}_0 & = &
a_{10}\widetilde{p}_0\left(\beta_0-\frac{1}{a}\beta_1\right) \\
{ } & = & a_{10}\widetilde{p}_0 \left(
\frac{a_{10}-a_{01}}{a_{00}a_{11}-a_{01}a_{10}} -
\frac{a_{10}b-a_{01}}{a(a_{00}a_{11}-a_{01}a_{10})}\right) \\
{ } & = &
\frac{a_{10}\widetilde{p}_0a_{00}\left(1-\frac{b}{a}\right)}{a_{00}a_{11}-a_{01}a_{10}} \\
{ } & \ne & 0.
\end{eqnarray*}
Similarly, we have
$\beta_0a_{11}\widetilde{p}_1-\frac{1}{b}\beta_1a_{11}\widetilde{p}_1\ne
0$. Therefore, $(q, (0,0))\in (\mathbb{R}^2)^2-\Sigma$. \hfill
$\Box$
\section{Proof of Theorem \ref{theorem 1}}\label{section 4}
\subsection{Proof of (1) of Theorem \ref{theorem 1}}\label{subsection 4.1}
Let $H_1:\R^2\to \R^2$ be the following diffeomorphism of the target
space.
\begin{eqnarray*}
{}&{}&H_1(X,Y)=\left(\frac{X}{a_{01}},\frac{Y}{a_{11}}\right)
\end{eqnarray*}
Set $r=\frac{a_{00}}{a_{01}}=\frac{a_{10}}{a_{11}}$. Then, the
composition of $G_{(p_0,p_1,A_1)}$ and $H_1$ can be expressed as
follows:
\begin{eqnarray*}
{}&{}&H_1\circ G_{(p_0,p_1,A_1)}(x,y)\\
{}&=&\left(\frac{a_{00}}{a_{01}}(x-p_{00})^2+(y-p_{01})^2,\frac{a_{10}}{a_{11}}(x-p_{10})^2+(y-p_{11})^2\right)\\
{}&=&\left(r(x-p_{00})^2+(y-p_{01})^2,r(x-p_{10})^2+(y-p_{11})^2\right).
\end{eqnarray*}
Let $H_2:\R^2\to \R^2$ be the following diffeomorphism of the source
space.
\begin{eqnarray*}
H_2(x,y)=\left(\frac{x}{\sqrt{|r|}}, y\right)
\end{eqnarray*}
The composition of $H_1\circ G_{(p_0,p_1,A_1)}$ and $H_2$ is as
follows:
\begin{eqnarray*}
{}&{}&H_1\circ G_{(p_1,p_2, A_1)}\circ H_2(x,y)\\
{}&=&\left(r\left(\frac{x}{\sqrt{|r|}}-p_{00}\right)^2+(y-p_{01})^2,r\left(\frac{x}{\sqrt{|r|}}-p_{10}\right)^2+(y-p_{11})^2\right).
\end{eqnarray*}
\par
Suppose that $r>0$.   Then, we have
\begin{eqnarray*}
{}&{}&H_1\circ G_{(p_0,p_1, A_1)}\circ H_2(x,y)\\
{}&=&((x-\sqrt{|r|}p_{00})^2+(y-p_{01})^2,(x-\sqrt{|r|}p_{10})^2+(y-p_{11})^2).
\end{eqnarray*}
Set $\widetilde{p}_0=(\sqrt{|r|}p_{00},p_{01})$ and
$\widetilde{p}_1=(\sqrt{|r|}p_{10},p_{11})$.    Then, we have
$H_1\circ G_{(p_0,p_1, A_1)}\circ
H_2=D_{(\widetilde{p}_0,\widetilde{p}_1)}$. Hence, by Proposition
\ref{proposition 1}, $G_{(p_0,p_1, A_1)}$ is
$\mathcal{A}$-equivalent to $\Phi_2$.
\par
Next, suppose that $r<0$.     Then, we have
\begin{eqnarray*}
{}&{}&H_1\circ G_{(p_0,p_1, A_1)}\circ H_2(x,y)\\
{}&=&(-(x-\sqrt{|r|}p_{00})^2+(y-p_{01})^2,-(x-\sqrt{|r|}p_{10})^2+(y-p_{11})^2).
\end{eqnarray*}
In this case, set $\widetilde{p}_0=(\sqrt{|r|}p_{00},p_{01})$ and
$\widetilde{p}_1=(\sqrt{|r|}p_{10},p_{11})$.   Then, we have
$H_1\circ G_{(p_0,p_1, A_1)}\circ H_2=L_{(\widetilde{p}_0,
\widetilde{p}_1)}$. Hence, by Proposition \ref{proposition 1},
$G_{(p_0,p_1, A_1)}$ is $\mathcal{A}$-equivalent to $\Phi_2$. \hfill
$\Box$
\subsection{Proof of (2) of Theorem \ref{theorem 1}}\label{subsection 4.2}
Since $G_{(p_0, p_1, A_2)}$ has a cusp point, it is not
$\mathcal{A}$-equivalent to $G_{(p_0, p_1, A_1)}$. Next, we show
that $G_{(p_0, p_1, A_2)}$ is proper. Suppose that $G_{(p_0, p_1,
A_2)}$ is not proper.   Then, by Proposition \ref{proposition 3},
$F_q$ is not proper, where $F_q$ is the mapping defined in
Proposition \ref{proposition 3}. Then, there exists a compact set
$C\subset \mathbb{R}^2$ such that $F_q^{-1}(C)$ is not bounded. Let
$\{\widetilde{p}_i\}_{i=1,2,\ldots}$ be a divergent sequence in
$F_q^{-1}(C)$. It follows that $F_q(\widetilde{p}_i)\in C$ for any
$i$ $(1\le i< \infty)$. Since $C$ is compact, taking a subsequence
if necessary, it follows that $\lim_{i\to
\infty}||F_q(\widetilde{p}_i)||<\infty$. On the other hand, by the
form of $F_q$, it also follows that $\lim_{i\to
\infty}||F_q(\widetilde{p}_i)||=\infty$. Hence, $G_{(p_0, p_1,
A_2)}$ is proper.
\par
By Proposition \ref{proposition 2} and Mather's characterization
theorem of proper stable mappings (\cite{mather5}), it is sufficient
to show that the restriction of $G_{(p_0, p_1, A_2)}$ to its
singular set $S(G_{(p_0, p_1, A_2)})$ is injective. In order to show
that  the restriction of $G_{(p_0, p_1, A_2)}$ to its singular set
$S(G_{(p_0, p_1, A_2)})$ is injective, by Proposition
\ref{proposition 3}, it is sufficient to show that the restriction
of $F_q$ to its singular set $S(F_q)$ is injective. Suppose the
restriction of $F_q$ to $S(F_q)$ is not injective. Then, since $F_q$
has the form $F_q(x,y)=((x-q_0)^2+(y-q_1)^2, ax^2+by^2)$, the point
$q=(q_0, q_1)$ must be inside the symmetry set of an ellipse defined
by $ax^2+by^2=c$ by some positive constant $c$. It is well-known
that the symmetry set of an ellipse defined by $ax^2+by^2=c$ must be
inside the union of two axes $\{(x, 0)\; |\; x\in \mathbb{R}\}$,
$\{(0, y)\; |\; y\in \mathbb{R}\}$ (for instance, see
\cite{brucegiblin}). However, this contradicts the fact that $(q,
(0,0))\in (\mathbb{R}^2)^2-\Sigma$.    Hence, the restriction of
$F_q$ to its singular set $S(F_q)$ is injective. \hfill $\Box$
\subsection{Proof of (3) of Theorem \ref{theorem 1}}\label{subsection 4.3}
By Proposition \ref{proposition 3}, $G_{(p_0,p_1,A_2)}$ is
$\mathcal{A}$-equivalent to
$F_{\widetilde{p}}:\mathbb{R}^2\to\mathbb{R}^2$ defined by
\begin{eqnarray*}
F_{\widetilde{p}}(x,y)=((x-\widetilde{p}_0)^2+(y-\widetilde{p}_1)^2,ax^2+by^2),
\end{eqnarray*}
where $a>0$, $b>0$, and
$(\widetilde{p},(0,0))\in(\mathbb{R}^2)^2-\Sigma $
$(\widetilde{p}=(\widetilde{p}_0,\widetilde{p}_1))$.

First, we show that there exist four non-zero real numbers $b_{00}$,
$b_{01}$, $b_{10}$, $b_{11}$ such that $F_{\widetilde{p}}$ is
$\mathcal{A}$-equivalent to the following mapping:
\begin{eqnarray*}
(x,y)\mapsto (x^2+b_{00}x+b_{01}y,y^2+b_{10}x+b_{11}y).
\end{eqnarray*}
Set
\begin{eqnarray*}
\widetilde{M}=\begin{pmatrix}
1 & a \\
1 & b
\end{pmatrix}^{-1}=
\frac{1}{a-b}\begin{pmatrix}
-b & a \\
1 & -1
\end{pmatrix}.
\end{eqnarray*}
The composition of the target diffeomorphism $(X,Y)\mapsto
(X,Y)\widetilde{M}$ and $F_{\widetilde{p}}$ has the following form:
\begin{eqnarray*}
(x,y)\mapsto
\left(x^2+\frac{2b\widetilde{p}_0}{a-b}x+\frac{2b\widetilde{p}_1}{a-b}y-\frac{b(\widetilde{p}_0^2+\widetilde{p}_1^2)}{a-b},
y^2-\frac{2a\widetilde{p}_0}{a-b}x-\frac{2a\widetilde{p}_1}{a-b}y+\frac{a(\widetilde{p}_0^2+\widetilde{p}_1^2)}{a-b}\right).
\end{eqnarray*}
By the target diffeomorphism to eliminate constant terms, we have
the following:
\[
(x,y)\mapsto
\left(x^2+\frac{2b\widetilde{p}_0}{a-b}x+\frac{2b\widetilde{p}_1}{a-b}y,
y^2-\frac{2a\widetilde{p}_0}{a-b}x-\frac{2a\widetilde{p}_1}{a-b}y\right).
\leqno{(4.1)}
\]
Since $ab\widetilde{p}_0\widetilde{p}_1\not=0$, all coefficients of
the variables $x$, $y$ are non-zero real numbers. Hence, there exist
four non-zero real numbers $b_{00}$, $b_{01}$, $b_{10}$, $b_{11}$
such that $F_{\widetilde{p}}$ is $\mathcal{A}$-equivalent to the
following mapping:
\[
(x,y)\mapsto (x^2+b_{00}x+b_{01}y,y^2+b_{10}x+b_{11}y).
\leqno{(4.2)}
\]
\par
Similarly, by Proposition \ref{proposition 3}, $G_{(q_0,q_1,B_2)}$
is $\mathcal{A}$-equivalent to
$F_{\widetilde{q}}:\mathbb{R}^2\to\mathbb{R}^2$ defined by
\begin{eqnarray*}
F_{\widetilde{q}}(x,y)=((x-\widetilde{q}_0)^2+(y-\widetilde{q}_1)^2,ax^2+by^2),
\end{eqnarray*}
where $a>0$, $b>0$, and
$(\widetilde{q},(0,0))\in(\mathbb{R}^2)^2-\Sigma $
$(\widetilde{q}=(\widetilde{q}_0,\widetilde{q}_1))$. Therefore,
there exist four non-zero real numbers $\widetilde{b}_{00}$,
$\widetilde{b}_{01}$, $\widetilde{b}_{10}$, $\widetilde{b}_{11}$
such that $F_{\widetilde{q}}$ is $\mathcal{A}$-equivalent to the
following mapping:
\[
(x,y)\mapsto
(x^2+\widetilde{b}_{00}x+\widetilde{b}_{01}y,y^2+\widetilde{b}_{10}x+\widetilde{b}_{11}y).
\leqno{(4.3)}
\]
\par
Finally, we show that the mapping defined by $(4.2)$ is
$\mathcal{A}$-equivalent to the mapping defined by $(4.3)$. For an
element $\beta =(\beta _{00},\beta _{01},\beta _{10},\beta _{11})\in
\mathbb{R}^4$ $(\beta _{00}\beta _{10}\not=0)$, we define the
following diffeomorphism of the target
\begin{eqnarray*}
H_\beta (X,Y)=\left(\frac{X}{\beta _{00}^2},\frac{Y}{\beta
_{10}^2}\right),
\end{eqnarray*}
and the following diffeomorphism of the source
\begin{eqnarray*}
h_\beta (x,y)=\left(\beta _{00}x+\beta _{01},\beta _{10}y+\beta
_{11}\right).
\end{eqnarray*}
By composing $H_\beta $, the mapping defined by $(4.1)$, and
$h_\beta $, we have the following:
\begin{eqnarray*}
(x,y)\mapsto \left(\frac{(\beta _{00}x+\beta _{01})^2+b_{00}(\beta _{00}x+\beta _{01})+b_{01}(\beta _{10}y+\beta _{11})}{\beta _{00}^2},\right.\\
\left. \frac{(\beta _{10}y+\beta _{11})^2+b_{10}(\beta _{00}x+\beta
_{01})+b_{11}(\beta _{10}y+\beta _{11})}{\beta _{10}^2}\right).
\end{eqnarray*}
By the target diffeomorphism to eliminate constant terms, we have
the following:
\begin{eqnarray*}
(x,y)\mapsto \left( x^2+\frac{2\beta _{01}+b_{00}}{\beta
_{00}}x+\frac{b_{01}\beta _{10}}{\beta _{00}^2}y,
y^2+\frac{b_{10}\beta _{00}}{\beta _{10}^2}x+\frac{2\beta
_{11}+b_{11}}{\beta _{10}}y\right),
\end{eqnarray*}
It is sufficient to show that there exists $\beta \in\mathbb{R}^4$
$(\beta _{00}\beta _{10}\not=0)$ satisfying
\[\frac{2\beta _{01}+b_{00}}{\beta _{00}}=\widetilde {b}_{00},\leqno{(4.4)}\]
\[\frac{b_{01}\beta _{10}}{\beta _{00}^2}=\widetilde {b}_{01},\leqno{(4.5)}\]
\[\frac{b_{10}\beta _{00}}{\beta _{10}^2}=\widetilde {b}_{10},\leqno{(4.6)}\]
\[\frac{2\beta _{11}+b_{11}}{\beta _{10}}=\widetilde {b}_{11}.\leqno{(4.7)}\]
By $(4.5)$ and $(4.6)$, we have
\begin{eqnarray*}
\beta _{00}=\left(\frac{b_{10}b_{01}^2}{\widetilde
{b}_{10}\widetilde {b}_{01}^2}\right)^{\frac{1}{3}}, \ \beta
_{10}=\left(\frac{b_{01}b_{10}^2}{\widetilde {b}_{01}\widetilde
{b}_{10}^2}\right)^{\frac{1}{3}}.
\end{eqnarray*}
Hence, it is clearly seen that there exist real numbers $\beta
_{01}$, $\beta _{11}$ satisfying $(4.4)$ and $(4.7)$. \hfill $\Box$
\par
\medskip
\noindent
{\bf Remark.}\quad
For a fixed matrix $A_2$ of rank 2 with non-zero entries,
Proposition 3 reduces to 2 the number of parameters of the mapping
$G_{(p_0,p_1,A_2)}$. Thus, for a fixed $p_0\in\mathbb R^2-\Sigma$,
we study the unfolding $\eta: \mathbb{R}^2\times \mathbb{R}^2\to
\mathbb{R}^2\times\mathbb{R}^2$ given by  $\eta(x,y,p_{10},
p_{11})=(G_{(p_0, p_1, A_2)}(x,y), p_{10}, p_{11})$. Based on Figure
\ref{figure 1}, it is expected that the unfolding has a $D_4^{+}$
singularity as a map germ.
\begin{definition}
{\rm Let $\Delta_{4}^+: \mathbb{R}^{4}\to \mathbb{R}^{4}$ denote the
following mapping:
\[
\Delta_{4}^+(x, y, u, v) =\left(x^2+u y, y^2+v x, u, v\right).
\]
When a map-germ $f: (\mathbb{R}^{4}, q)\to (\mathbb{R}^{4}, f(q))$
is $\mathcal{A}$-equivalent to $\Delta_{4}^+: (\mathbb{R}^{4}, 0)\to
(\mathbb{R}^{4},0)$, the point $q\in \mathbb{R}^{4}$ is said to be a
{\it $D_4^+$ singularity of $f$}. }
\end{definition}
\noindent The expectation can be easily realized by the first half
of the proof of (3) of Theorem \ref{theorem 1}. This is because
(4.1) is $\mathcal{A}$-equivalent to
\[
(x,y)\mapsto \left(x^2+\frac{2b\widetilde{p}_1}{a-b}y,
y^2-\frac{2a\widetilde{p}_0}{a-b}x\right)
\]
and thus the mapping
\[
(x, y, \widetilde{p}_0, \widetilde{p}_1) \mapsto (F_q(x,y),
\widetilde{p}_0, \widetilde{p}_1)
\]
is $\mathcal{A}$-equivalent to
 \[
(x,y, \widetilde{p}_0, \widetilde{p}_1) \mapsto
\left(x^2+\widetilde{p}_1 y, y^2+\widetilde{p}_0 x,
\widetilde{p}_{0}, \widetilde{p}_{1}\right).
\]
On the other hand, by the proof of Proposition \ref{proposition 3},
it follows that
\begin{eqnarray*}
\widetilde{p}_0 & = &
\beta_0a_{10}(p_{10}-{p}_{00})- \frac{1}{a}\beta_1a_{10}(p_{10}-{p}_{00}),\\
\widetilde{p}_1 & = &
\beta_0a_{11}(p_{11}-{p}_{01})-\frac{1}{b}\beta_1a_{11}(p_{11}-{p}_{01}).
\end{eqnarray*}
Hence, we have the following:
\begin{proposition}\label{proposition 4}
Let $p_0$ be a given point of $\mathbb{R}^2$ and
let $A_2$ be a $2\times 2$ matrix of rank $2$ with non-zero entries.
Then,
$\eta: \mathbb{R}^2\times \mathbb{R}^2\to
\mathbb{R}^2\times\mathbb{R}^2$ is $\mathcal{A}$-equivalent to
$\Delta_4^+$, and thus the point $(p_0, p_0)\in \mathbb{R}^4$ is a
$D_4^+$ singularity of $\eta$.
\end{proposition}
\section{Geometric interpretation of Propositions \ref{proposition 2} and \ref{proposition 3}}\label{section 5}
Let $A_2$ be a $2\times 2$ matrix of rank 2. Then the generalized
distance-squared mapping $G_{(p_0,p_1,A_2)}=(G_1,G_2)$ determines
two foliations in the plane by conics $$\mathcal
C_1(c_1)=\{(x,y)\,|\,G_1(x,y)=c_1\} \, \text{and}\,\,\mathcal
C(c_2)=\{(x,y)\,|\,G_2(x,y)=c_2\},$$ where $c_1, c_2 \in \mathbb R.$
The two foliations are generically transverse, and they are tangent
at the singularities of $G_{(p_0,p_1,A_2)}$. It follows from
Proposition 3  that $G_{(p_0,p_1,A_2)}$ is equivalent to
$F_q(x,y)=\left((x-q_0)^2+(y-q_1)^2, ax^2+by^2\right),$ for some
point $q=(q_0, q_1)$ and positive constants $a$ and $b.$  On the
other hand, the singularities of $F_q$ correspond to the
singularities of the family of distance squared functions
$D_q(x,y)=(x-q_0)^2+(y-q_1)^2$ on the family of ellipses
$\{ax^2+by^2=c_2\}$.

The geometric characterizations of the singularities of the family
of distance-squared functions on a unit speed plane curve $\alpha: I
\to \mathbb R^2,\, \alpha(s)=(x(s), y(s))$ are well known, and we
recall them in the following proposition. As usual, we write $\{{\bf
t}(s), {\bf n}(s)\}$ for the Frenet frame  and $k(s)$ the curvature
of $\alpha$  at the point $s.$

\begin{proposition}[\cite{brucegiblin}, page 29]
\label{proposition 5} Let $f_q(s)=D_q(\alpha(s)). $ Then, the
following hold
\begin{enumerate}
\item $f_q'(s_0)=0$ if and only if there exists $\lambda\in \mathbb{R}$ such that
$ q=\alpha(s_0)+\lambda {\bf n}(s_0),$
\item $f_q'(s_0)=f_q''(s_0)=0$ if and only if $k(s_0)\neq 0$ and
$q=\alpha(s_0)+ \frac{{\bf n}(s_0)}{k(s_0)},$
\item $f_q'(s_0)=f_q''(s_0)=f_q'''(s_0)=0$ if and only if
$k(s_0)\neq 0,$  $q=\alpha(s_0)+ \frac{{\bf n}(s_0)}{k(s_0)}$ and
$k'(s_0)=0.$
\end{enumerate}
\end{proposition}

In the previous setting the point $q=\alpha(s_0)+ \frac{{\bf
n}(s_0)}{k(s_0)}$ is the center of curvature of the curve $\alpha$
at $s=s_0$. The circle centered at $q$ passing through $\alpha(s_0)$
is the osculating circle of the curve at this point.

These characterizations and the proof of Proposition 2 yield the
following

\begin{proposition}
\label{proposition 6}
\begin{enumerate}
\item A point $\rho=(x_0,y_0)\in\mathbb R^2$ is a fold point of $F_q$ if and
only if the circle $\{(x-q_0)^2+(y-q_1)^2=(x_0-q_0)^2+(y_0-q_1)^2\}$
and the ellipse $ax^2+by^2=ax_0^2+by_0^2$ are tangent at $\rho$ and
the center of curvature of the ellipse at $\rho$ is different from
$(q_0,q_1)$.

\item A point $\rho=(x_0,y_0)\in\mathbb R^2$ is a cusp point of $F_q$ if and
only if the circle $\{(x-q_0)^2+(y-q_1)^2=(x_0-q_0)^2+(y_0-q_1)^2\}$
and the ellipse $ax^2+by^2=ax_0^2+by_0^2$ are tangent at $\rho$ and
the center of curvature of the ellipse is $(q_0,q_1)$.
\end{enumerate}
\end{proposition}

\begin{corollary}
Given any two families of conics $\mathcal C_1(c_1)$ and $\mathcal
C_2(c_2)$ as above defining a mapping $G_{(p_0,p_1,A)},$
$\text{rank}\, A=2,$ there exist  uniquely determined constants
$c_1$ and $c_2$ and a  point $\rho=(x_0,y_0)$ such that
$\{\rho\}=\mathcal C_1(c_1)\cap \mathcal C_2(c_2)$ is the unique
point of tangency between the two curves  at which the osculating
circles coincide.
\end{corollary}
\section*{Acknowledgements}
S.~Ichiki and T.~Nishimura are partially supported by JSPS-CAPES
under the JAPAN-BRAZIL research cooperative program. R.~Oset Sinha
is partially supported by FAPESP grant no. 2013/02381-1 and DGCYT
and FEDER grant no. MTM2012-33073. M. A. S.~Ruas is partially
supported by CNPq grant no. 305651/2011-0.


\end{document}